\documentclass[pdflatex,sn-mathphys-num]{sn-jnl}
\usepackage{cmap}
\usepackage[USenglish]{babel}
\usepackage[T1]{fontenc}

\usepackage[utf8]{inputenc}				
\usepackage{microtype}

\usepackage{graphicx}
\usepackage{multirow}
\usepackage{subcaption}

\usepackage{amsmath,amssymb,amsfonts}
\usepackage{amsthm}

\usepackage{xcolor}

\theoremstyle{thmstyleone}
\newtheorem{theorem}{Theorem}

\theoremstyle{thmstylethree}

\theoremstyle{thmstyletwo}

\usepackage{algorithm}
\usepackage{algpseudocode}
\usepackage{booktabs}

\begin{document}

\title[Forward Error-Oriented Iterative Refinement for Eigenvectors]{Forward Error-Oriented Iterative Refinement for Eigenvectors of a Real Symmetric Matrix}

\author*[1]{\fnm{Takeshi} \sur{Terao}}
\author[2]{\fnm{Katsuhisa} \sur{Ozaki}}

\affil*[1]{Research Institute for Science and Engineering, Waseda University, 3-4-1 Okubo, Shinjuku-ku, Tokyo 169-8555, Japan}
\affil[2]{Department of Mathematical Sciences, Shibaura Institute of Technology, 307 Fukasaku, Minuma-ku, Saitama, 337-8570, Saitama, Japan}

\abstract{
In this paper, we discuss numerical methods for the eigenvalue decomposition of real symmetric matrices.
While many existing methods can compute approximate eigenpairs with sufficiently small backward errors, the magnitude of the resulting forward errors is often unknown.
Consequently, when high-precision numerical solutions are required, the computational cost tends to increase significantly because backward errors must be reduced to an excessive degree.
To address this issue, we propose an efficient approximation algorithm that aims to achieve a prescribed forward error, together with a high-accuracy numerical algorithm based on the Ozaki scheme---an emulation technique for matrix multiplication---adapted to this problem.
Since the proposed method is not primarily focused on reducing backward errors, the computational cost can be significantly reduced.
Finally, we present numerical experiments to evaluate the efficiency of the proposed method.
  }

\maketitle

\section{Introduction}

In scientific computing and machine learning, linear algebra computations—such as solving systems of linear equations and eigenvalue problems—often account for the majority of the overall computational cost.
In response, numerous high-performance algorithms have been developed to address these challenges.
Although numerical computations are generally efficient, their finite precision inevitably introduces rounding errors.
To evaluate the accuracy and reliability of such computations, it is essential to understand both forward and backward errors and to consider their respective characteristics.

We consider the eigenvalue problem of a real symmetric matrix $A \in \mathbb{R}^{n \times n}$, defined by
\begin{align}
    A x_{(i)} = \lambda_i x_{(i)}, \quad \|x_{(i)}\|_2 = 1, \quad i = 1, 2, \dots, n,
\end{align}
where $\lambda_i \in \mathbb{R}$ and $x_{(i)} \in \mathbb{R}^n$ denote the eigenvalue and the corresponding eigenvector, respectively.
The eigenvalue matrix $D \in \mathbb{R}^{n \times n}$ and the eigenvector matrix $X \in \mathbb{R}^{n \times n}$ are defined by
\begin{align}
    D = \mathrm{diag}(\lambda_1, \lambda_2, \dots, \lambda_n), \quad
    X = \left(x_{(1)}, x_{(2)}, \dots, x_{(n)}\right).
\end{align}
For this problem, the following relationships hold:
\begin{align}
    AX=XD,\quad X^TX=I,
\end{align}
where $I$ denotes the identity matrix of appropriate size.
The approximate eigenvalue matrix and eigenvector matrix are denoted by $\widehat D$ and $\widehat X$, respectively.

In many cases, the accuracy of the eigenvalue decomposition is assessed using backward errors, such as orthogonality $\|I-\widehat X^T\widehat X\|$ and diagonalization $\|\widehat X^T A \widehat X - \widehat D\|$.  
Most numerical libraries ensure that these backward errors are sufficiently small.  
However, in practice, it is often more important to obtain approximate solutions with sufficiently small forward errors, such as $\|X - \widehat X\|$.

\begin{figure}[htbp]
    \centering
    \includegraphics[width=0.5\textwidth]{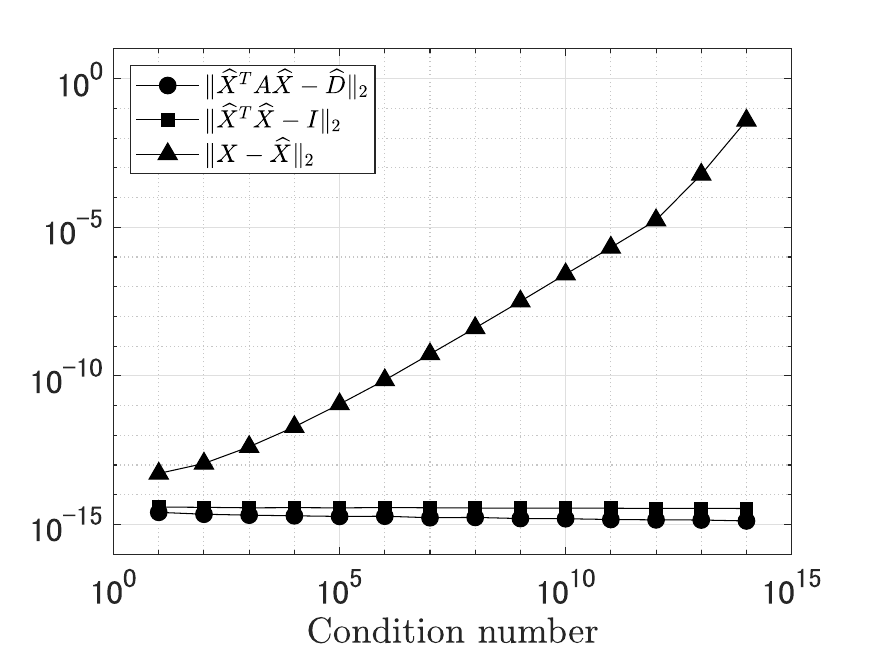}
    \caption{Comparison of forward error $\|X-\widehat X\|_2$, and backward errors $\|\widehat X^T\widehat X-I\|_2$ and $\|\widehat X^TA\widehat X-\widehat D\|_2$ for an eigenvector matrix ($n=100$). The condition number is $\|A\|_2\cdot\|A^{-1}\|_2$.}
    \label{fig1}
\end{figure}

The difficulty lies in the fact that, even when the backward error is sufficiently small, the forward error may still be large.
We conducted the following test using MATLAB:
\begin{verbatim}
    A = gallery('randsvd',n,-cnd);
    [X,D] = eig(A);
\end{verbatim}
Here, $A$ is an $n \times n$ real symmetric positive definite matrix with a condition number specified by \texttt{cnd}.  
The eigenvalues of $A$ are distributed geometrically.  
The variables \texttt{X} and \texttt{D} are the approximate eigenvector matrix and the approximate eigenvalue matrix, respectively.  
Results comparing the forward and backward errors are shown in Figure~\ref{fig1}.  
The reference solution $X$ is computed using a pseudo-quadruple precision numerical computation library~\cite{advanpix2006multiprecision}.

It can be observed from Figure~\ref{fig1} that small backward errors are obtained even for ill-conditioned matrices.
However, the forward error grows as the condition number increases.  
In practice, since the true solution is typically unavailable, the magnitude of the forward error often remains unknown.

The iterative refinement method proposed by Ogita and Aishima~\cite{ogita2018iterative,ogita2019iterative} has recently attracted considerable attention, particularly in combination with matrix multiplication based on error-free transformations~\cite{ozaki2012error,ozaki2013generalization,ozaki2024extension}.  
Iterative refinement methods are now recognized as high-performance mixed-precision numerical algorithms~\cite{abdelfattah2021survey,higham2022mixed}, and matrix multiplication using error-free transformations is being actively developed as an efficient computational technique for current and emerging architectures~\cite{ootomo2024dgemm,uchino2024performance,dongarra2024hardware,schwarz2025guaranteed}.  
This combination has also gained attention as a high-performance and high-accuracy eigenvalue solver~\cite{uchino2024high}.

\begin{figure}[htbp]
    \centering
    \begin{subfigure}{0.49\textwidth}
        \centering
        \includegraphics[width=\textwidth]{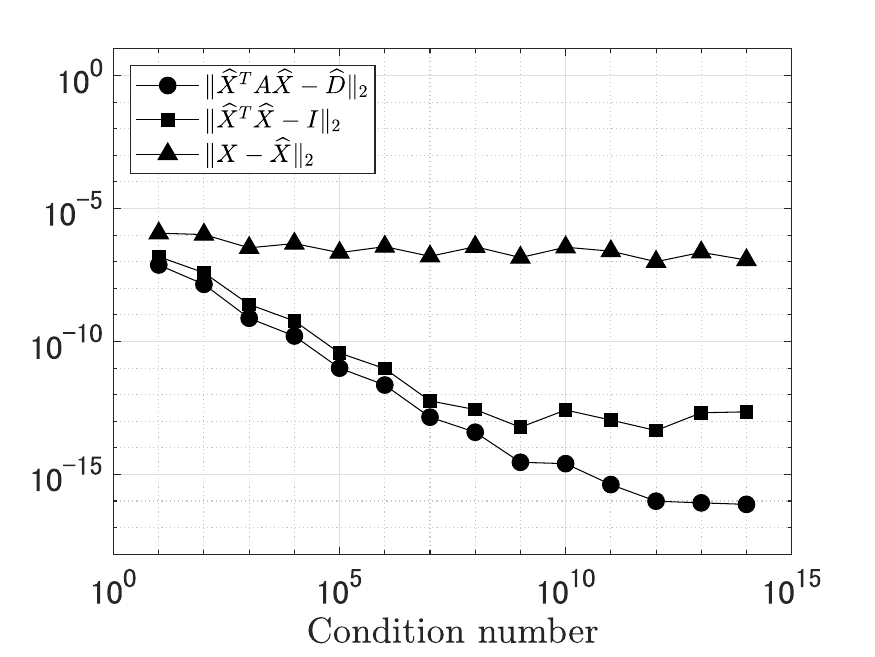} 
        \subcaption{$\delta=10^{-6}$}\label{fig2}
    \end{subfigure}
    \hfill
    \begin{subfigure}{0.49\textwidth}
        \centering
        \includegraphics[width=\textwidth]{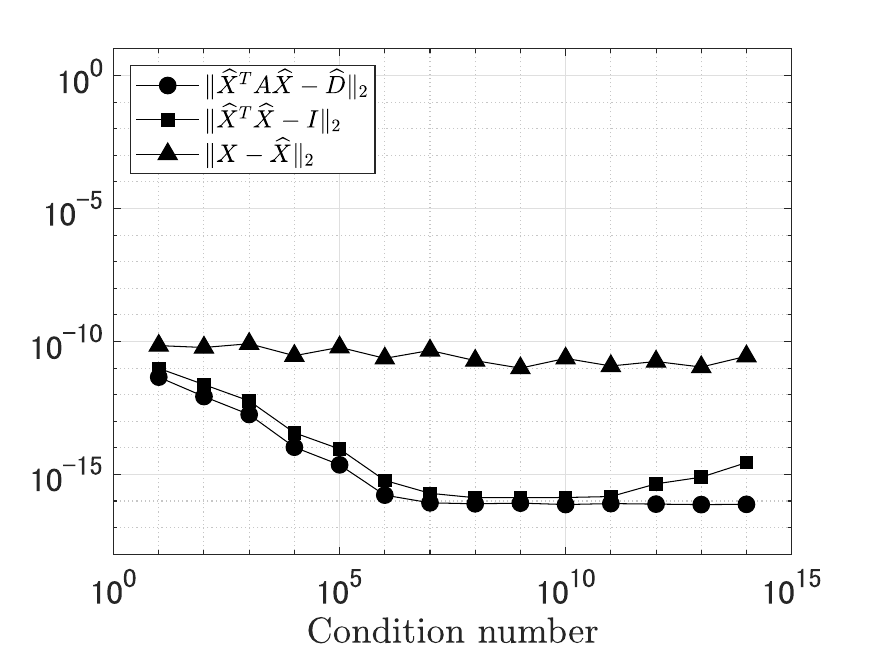} 
        \subcaption{$\delta=10^{-10}$}\label{fig3}
    \end{subfigure}
    \caption{Goal results of proposed method $(n=100)$. Here, $\delta$ is the specifiable scalar, and the condition number is $\|A\|_2\cdot\|A^{-1}\|_2$.}
\end{figure}

The present focuses on forward error control and seeks to reduce the computational cost of high-precision matrix multiplications based on the Ogita--Aishima method.  
In addition, we adapt the Ozaki scheme to effectively integrate it into the proposed method.  
As a result, the proposed method enables efficient computation of eigenvectors that approximately satisfy a prescribed forward error tolerance:
given a constant $\delta$, the method computes $\widehat{X}$ such that $\|X - \widehat{X}\| \lesssim \delta$.  
Figures~\ref{fig2} and \ref{fig3} present computational results for $\delta = 10^{-6}$ and $\delta = 10^{-10}$, respectively.  
If $\|X - \widehat{X}\| > \delta$, then the approximate solution can be refined to meet the specified tolerance.

The remainder of this paper is organized as follows.  
Section~2 introduces the notation used throughout the paper and reviews the iterative refinement method for eigenvalue decomposition, together with the high-precision matrix multiplications that form the foundation of the proposed approach.  
Section~3 describes a refined iterative method with reduced computational cost, building on the techniques described in Section~2.  
Section~4 describes numerical experiments that demonstrate the effectiveness of the proposed method.  
Finally, Section~5 summarizes the main features of the method and outlines directions for future work.

\section{Previous work}

\subsection{Notation}

We first give some notation used in this paper.
The floor function, denoted by $\lfloor x \rfloor$, returns the greatest integer less than or equal to $x$, that is, 
\[
\lfloor x \rfloor := \max\{ n \in \mathbb{Z} \mid n \leq x \}.
\]
Similarly, the ceiling function $\lceil x \rceil$ is defined as the smallest integer that is greater than or equal to $x$, that is,
\[
\lceil x \rceil := \min\{ n \in \mathbb{Z} \mid n \geq x \}.
\]
Unless otherwise specified, the 2-norm is used for both matrices and vectors.

Let $\mathbb{F}$ denote the set of floating-point numbers, and let $u$ be the unit roundoff of $\mathbb{F}$ (e.g., $u=2^{-53}$ for binary64).
The $(i,j)$ element of a matrix $A$ is denoted by $a_{ij}$, and the $i$-th element of a vector $v$ is denoted by $v_i$.
For a real number $\alpha$, we define
\begin{align}
    \mathrm{ufp}(\alpha):=
    \begin{cases}
        2^{\lfloor \log_2|\alpha| \rfloor}, & \alpha \neq 0, \\
        0, & \alpha = 0.
    \end{cases}
\end{align}
The meaning of $\mathrm{ufp}(\alpha)$ is a unit in the first place of $\alpha$.
$\mathrm{fl}(\cdot)$ denotes the result of evaluating all operations within the parentheses in floating-point arithmetic under the round-to-nearest mode.

\subsection{Pairwise arithmetic}

The proposed method requires retaining the computation results with a higher precision than that of standard floating-point numbers.
Therefore, this section introduces the necessary algorithms for the computation.

For $a, b \in \mathbb{F}$, we can obtain $x$ and $y$ such that
\[
x = \mathrm{fl}(a+b), \quad y = a + b - x.
\]
The above $y$ is obtained by 
\[
  y = \mathrm{fl}(( a - (x-z) ) + (b-z)), 
\]
where $z = \mathrm{fl}(x - a)$.
This algorithm was proposed by Knuth~\cite{knuth1997seminumerical}.
Denote this algorithm by $[x,y] = \mathtt{TwoSum}(a,b)$, which follows the MATLAB notation.

Let $b_h, b_\ell \in \mathbb{F}$ with the assumption $|b_h| \gg |b_\ell|$.
We obtain an approximation of $a + b_h + b_\ell$ by $c_h + c_\ell, \ c_h, c_\ell \in \mathbb{F}$ such that
\[
[c_h,t] = \mathtt{TwoSum}(a,b_h), \quad c_\ell = \mathrm{fl}(b_\ell + t).
\]
This is a pair arithmetic proposed in~\cite{lange2020faithfully} 
that omits the lower part of another double-word number.

\subsection{Error-free transformation for matrix multiplication}

Here, we introduce high-precision matrix computation methods using error-free transformations~\cite{ozaki2012error,ozaki2013generalization,ozaki2024extension}, which are collectively called the Ozaki scheme.
For the matrices $\underline A^{(1)}:=A$ and $\underline{\widehat X}^{(1)}:={\widehat X}$, let
\begin{align}
    v_i^{(r)}:=\max_j|\underline a_{ij}^{(r)}|,\quad w_j^{(r)}:=\max_i|\underline{\widehat x}_{ij}^{(r)}|.
\end{align}
For $\alpha, \beta \in \mathbb{N}$, define
\begin{align}
    \sigma_i^{(r)}&:=
    \begin{cases}
        0, & v_i^{(r)}=0, \\
        0.75\cdot 2^{\lceil \log_2 v_i^{(r)} \rceil}\cdot 2^{\alpha}, & \text{otherwise},\\
    \end{cases}\label{eq:sigma} \\
    \tau_j^{(s)}&:=
    \begin{cases}
        0, & w_j^{(s)}=0, \\
        0.75\cdot 2^{\lceil \log_2 w_j^{(s)} \rceil}\cdot 2^{\beta}, & \text{otherwise}
    \end{cases}\label{eq:tau}
\end{align}
from \cite{minamihata2016improved}.
Using these definitions, for $A^{(r)}$, we compute
\begin{align}
    a_{ij}^{(r)}=\mathrm{fl}\left(\left(\sigma_i^{(r)}+\underline a_{ij}^{(r)}\right)-\sigma_i^{(r)}\right),\quad \underline a_{ij}^{(r+1)}=\mathrm{fl}\left(\underline a_{ij}^{(r)}-a_{ij}^{(r)}\right) \label{eq:a}
\end{align}
for $r=1,2,\dots$\ .
Similarly, for $\widehat X$, we compute
\begin{align}
    \widehat x_{ij}^{(s)}=\mathrm{fl}\left(\left(\tau_j^{(s)}+\underline{\widehat x}_{ij}^{(s)}\right)-\tau_j^{(s)}\right),\quad \underline{\widehat x}_{ij}^{(s+1)}=\mathrm{fl}\left(\underline{\widehat x}_{ij}^{(s)}-\widehat x_{ij}^{(s)}\right)\label{eq:x}
\end{align}
for $s=1,2,\dots$\ .
For these terms, it has been shown that
\begin{align}
    |\underline a_{ij}^{(r+1)}| \leq u \cdot \mathrm{ufp}\left(\sigma_i^{(r)}\right),\quad |\underline{\widehat x}_{ij}^{(s+1)}| \leq u \cdot \mathrm{ufp}\left(\tau_j^{(r)}\right)\label{eq:minamihata}.
\end{align}
Additionally, if
\begin{align}
    \alpha+\beta \geq -\log_2u+\log_2n,\label{eq:alp_beta}
\end{align}
then
\begin{align}
    A^{(r)}\widehat X^{(s)}=\mathrm{fl}\left(A^{(r)}\widehat X^{(s)}\right),\quad r=1,2,\dots, n_A-1, \quad s=1,2,\dots, n_X-1
\end{align}
holds.
Therefore, if we set sufficiently large $n_A$ and $n_X$, then
\begin{align}
\underline{A}^{n_A+1} = {\bf O}, \quad \underline{\widehat{X}}^{n_X+1} = {\bf O}
\end{align}
and 
\begin{align}
    A\widehat X = \sum_{r=1}^{n_A}\sum_{s=1}^{n_X} A^{(r)}\widehat X^{(s)} = \sum_{r=1}^{n_A}\sum_{s=1}^{n_X}\mathrm{fl}\left(A^{(r)}\widehat X^{(s)}\right)
\end{align}
are valid.
Here, \({\bf O}\) denotes a zero matrix of appropriate size.
As an efficient method, when $k := n_A = n_X$, matrix multiplication using Ozaki-scheme can be performed with $\frac{1}{2}k(k+1)$ matrix multiplications and $\frac{1}{2}k(k+1)-1$ matrix additions~\cite{ozaki2013generalization}.

Some typical forms that involve 3, 6, and 10 matrix multiplications are as follows:
\begin{align}
AB & = A^{(1)} B^{(1)} + A^{(1)} B^{(2)} + \underline{A}^{(2)} B \label{eq:mat3}\\  
 & = A^{(1)} B^{(1)} + A^{(1)} B^{(2)} + A^{(2)} B^{(1)} + A^{(2)} B^{(2)} + \underline{A}^{(3)} B + (A^{(1)} + A^{(2)}) \underline{B}^{(3)} \label{eq:mat6}\\   
 & = A^{(1)} B^{(1)} + A^{(1)} B^{(2)} + A^{(2)} B^{(1)} + A^{(1)} B^{(3)} + A^{(2)} B^{(2)} + A^{(3)} B^{(1)} \nonumber \\
 & \quad + \left(A - A^{(4)}\right) \underline{B}^{(4)} + \left( A^{(2)} + A^{(3)} \right) B^{(3)} + A^{(3)} B^{(2)} +  \underline{A}^{(4)} B. \label{eq:mat10}  
\end{align}
See~\cite{ozaki2024extension} for how to set $\alpha$ and $\beta$ ($\alpha = \beta$) for each form.

\subsection{Iterative refinement for symmetric eigenvalue decomposition}

An efficient iterative refinement method for the eigenvalue decomposition of real symmetric matrices was proposed by Ogita and Aishima~\cite{ogita2018iterative,ogita2019iterative}, and further developed for error analysis and cost reduction~\cite{Shirama2019,uchino2022}.

Let $A \in \mathbb{R}^{n \times n}$ be a real symmetric matrix with distinct eigenvalues, and let $X, \widehat X \in \mathbb{R}^{n \times n}$ be the exact and approximate eigenvector matrices, respectively.
Supposing $X = \widehat X(I+E)$ for $\|E\| \leq 1/\sqrt{2}$, and using the orthogonality and diagonalization conditions $X^T X = I$, $X^T A X = D$, one obtains the relations
\begin{align}
    \widehat X^T \widehat X &= (I+E)^{-T}(I+E)^{-1},\\
    \widehat X^T A \widehat X &= (I+E)^{-T} D (I+E)^{-1}.
\end{align}
Neglecting higher-order terms in the Neumann series expansion, the correction matrix $E$ is approximated by solving
\begin{align}
    \widetilde E + \widetilde E^T &= I - \widehat X^T \widehat X =: R,\\
    \widetilde D - \widetilde D \widetilde E - \widetilde E^T \widetilde D &= \widehat X^T A \widehat X =: S,
\end{align}
where $\widetilde D = \mathrm{diag}(\widetilde \lambda_1, \dots, \widetilde \lambda_n)$ is computed from the Rayleigh quotients using $\widehat{X}$ and $A$.
Solving these equations element-wise yields
\begin{align}
    \widetilde e_{ij} =
    \begin{cases}
        \dfrac{r_{ii}}{2}, & i = j, \\
        \dfrac{s_{ij} + \widetilde \lambda_j r_{ij}}{\widetilde \lambda_j - \widetilde \lambda_i}, & \text{otherwise}.
    \end{cases}
\end{align}
The updated approximation is given by $\widehat X \leftarrow \widehat X + \widehat X \widetilde E$.
For $W=\widehat X^T(A\widehat X-\widehat X\widetilde D)$, it holds that
\begin{align}
    w_{ij}=s_{ij}+\widetilde \lambda_j r_{ij}
\end{align}
for $i\not=j$.
\begin{algorithm}
\caption{Refinement algorithm based on Ogita-Aishima method (one iteration).}
\begin{algorithmic}[1]
    \Require $A = A^T, \widehat X \in \mathbb{R}^{n \times n}$
    \Ensure $\widehat X \in \mathbb{R}^{n \times n}$
    \State $\displaystyle V = \mathrm{accmul}(A,\widehat X)$\hfill \text{\% Accurate matrix multiplication, for example, see Section 2.3}
    \State $r_i = 1-\widehat x_i^T \widehat x_i$ for $i=1,2,\dots,n$
    \State $\widehat \lambda_i = \widehat x_i^T \widehat v_i / (1-r_i)$ for $i=1,2,\dots,n$
    \State $W = \widehat X^T (V - \widehat X\widehat D)$, where $\widehat D = \mathrm{diag}(\widehat \lambda_1, \dots, \widehat \lambda_n)$
    \State $\widetilde e_{ij} =
    \begin{cases}
        \frac{w_{ij}}{\widetilde \lambda_j - \widetilde \lambda_i}, & i \neq j, \\
         r_{i} / 2, & i = j
    \end{cases}$
    \State $\widehat X \leftarrow \mathrm{accsum}(\widehat X, \widehat X\widetilde E)$\hfill \text{\% Accurate summation, for example, Section 2.2}
\end{algorithmic}

\end{algorithm}

The following theorem provides an upper bound on the error of the computed correction matrix $\widetilde{E}$ obtained from Algorithm~1~\cite{Shirama2019}.
This theorem guarantees the validity of the update $\widehat X \leftarrow \widehat X + \widehat X \widetilde E$, under the assumption that the initial approximation $\widehat X$ is sufficiently accurate.

\begin{theorem}\label{thm:yamamoto}
Let $A \in \mathbb{R}^{n \times n}$ be a real symmetric matrix with all distinct eigenvalues. Apply Algorithm~1 to $A$ and $\widehat X \in \mathbb{R}^{n \times n}$.  
Let $\epsilon := \|E\|$, and let $e_j$ be the $j$-th column of $E$. If
\begin{align}
    \epsilon < \min\left(\frac{\min_{i \neq j}|\lambda_j - \lambda_i|}{10\sqrt{n}\|A\|}, \frac{1}{100}\right),
\end{align}
then
\begin{align}
    \|E - \widetilde E\| \leq \frac{6.4\sqrt{n}\|A\|\epsilon^2}{\min_{i \neq j}|\lambda_j - \lambda_i| - 14\|A\|\epsilon}.
\end{align}
\end{theorem}

This theorem implies that, when the eigenvalues of $A$ are well-separated (i.e., there are no clusters or multiplicities), 
the correction matrix $\widetilde{E}$ approximates the true error $E$ with an accuracy of $\mathcal{O}(\epsilon^2)$. Consequently, the iterative process exhibits quadratic convergence, provided that the initial error is sufficiently small.

\section{Main result}

\subsection{Outline}
We describe a method for computing an approximation $\widetilde X$ of the eigenvector matrix $X$ for a symmetric matrix $A \in \mathbb{R}^{n \times n}$ such that
\begin{align}
    \|X - \widetilde X\| \approx \delta > u,
\end{align}
where $\delta$ is a user-specified constant.

In the original iterative refinement, a high-precision computation of the matrix product $A\widehat{X}$ is required (see, for example, Section~2.2).  
Assuming that we use $\widehat{X} = \widehat{X}^{(1)} + \widehat{X}^{(2)}$ for accurate computations, the proposed method neglects $\widehat{X}^{(2)}$ in the decomposition.  
In other words, the proposed method was constructed to compute the matrix product $A\widehat{X}^{(1)}$ with high precision.  
Note that $A\widehat{X}^{(1)} \neq A\widehat{X}$.

When $A^{(r)}$ and $\widehat{X}^{(s)}$ are defined as in \eqref{eq:a} and \eqref{eq:x}, respectively, the following holds:
\begin{align}
    \sum_{r=1}^{n_A} A^{(r)} \widehat{X}^{(1)} = \sum_{r=1}^{n_A} \mathrm{fl}(A^{(r)} \widehat{X}^{(1)}),
\end{align}
which implies that the number of standard matrix multiplications required to compute the high-precision matrix product is $n_A$.

A key idea of the proposed method lies in how the constant $\beta$ in \eqref{eq:tau}, used to compute $\widehat{X}^{(1)}$, is determined.  
When $\beta$ is small, $\alpha$ increases accordingly, and a larger $n_A$ is required.
As a result, the computational cost of the high-precision matrix product increases.  
On the other hand, if $\beta$ is too large, it can become impossible to compute an $\widetilde{X}$ that satisfies $\|X - \widetilde{X}\| \approx \delta$.  
To compute a suitable $\beta$, we define constants and assumptions in the next section and perform an error analysis using them.

\subsection{Definitions and assumptions}

The matrices $X$ and $\widehat X$ denote the eigenvector matrix and its approximation of $A=A^T$, respectively.
Suppose that the correction matrix $E$ and its approximation $\widetilde E$ are such that
\begin{align}
    X = \widehat X^{(1)}(I+E),\quad X' := \widehat X^{(1)}(I+\widetilde E),\quad \widehat X = \widehat X^{(1)} + \widehat X^{(2)}.\label{eq:E'}
\end{align}
Moreover, we consider eigenvalue problems without clustered eigenvalues.
In addition, we assume that the condition in Theorem~\ref{thm:yamamoto}
\begin{align}
    \epsilon := \|E\| <  \min\left(\frac{\min_{i \neq j}|\lambda_j - \lambda_i|}{10\sqrt{n}\|A\|}, \frac{1}{100}\right)\label{eq:assume1}
\end{align}
is satisfied for the matrix $E$.
Furthermore, we assume that there exists $\xi$ such that
\begin{align}
    \|E - \widetilde E\| &= \epsilon^2\cdot \frac{1}{\xi} \cdot \frac{\max_k|\widehat \lambda_k|}{\min_{i \neq j}|\widehat \lambda_j - \widehat \lambda_i|}, \quad \min_{i \neq j}|\widehat \lambda_j - \widehat \lambda_i| \not = 0,\quad \xi\not=0.\label{eq:assume2}
\end{align}

We now explain the necessity of the assumption~\eqref{eq:assume1}.
The proposed method relies on the inequality~\eqref{eq:assume2}.
Therefore, if the inequality changes, the splitting method must be reconsidered.
If the assumption~\eqref{eq:assume1} is not satisfied, a different inequality must be used~\cite{ogita2018iterative,ogita2019iterative,Shirama2019,uchino2024acceleration}.
Although the techniques of the proposed method may still be applicable to such cases, the discussion would become complex and lengthy, so we omit it from this paper.

\subsection{Matrix splitting method}

Here, we describe a method for computing the required constant $\beta$ for the slice $\widehat{X} = \widehat{X}^{(1)} + \widehat{X}^{(2)}$ of the matrix $\widehat{X}$ under reasonable assumptions. 

First, we describe the conditions that $\widehat X^{(2)}$ must satisfy.
Assuming $\|(\widehat X^{(1)})^T\widehat X^{(1)}\| \approx 1$, we obtain
\begin{align}
    \delta &\approx \|X - X'\| = \|X - \widehat X^{(1)}(I+\widetilde E)\| \\
    &= \|X - \widehat X^{(1)}(I+E) + \widehat X^{(1)}(E - \widetilde E)\| =\|\widehat X^{(1)}(E - \widetilde E)\|\\
    &\approx \|E - \widetilde E\|\label{aprx_E}.
\end{align}
From this, consider the following conditions 
\begin{align}
    \begin{cases}
        \|X-X'\|\leq \delta\\
        \|\widehat X^{(2)}\|\leq \epsilon
    \end{cases}.
\end{align}
Thus, from \eqref{eq:assume2} and \eqref{aprx_E}, we require
\begin{align}
    \epsilon^2 \cdot\frac{1}{\xi}\cdot \frac{\max_k|\widehat \lambda_k|}{\min_{i \neq j}|\widehat \lambda_j - \widehat \lambda_i|} \leq \delta\quad \Longrightarrow \quad 
    \epsilon \leq \sqrt{\delta\xi\cdot\frac{\min_{i \neq j}|\widehat \lambda_j - \widehat \lambda_i|}{\max_k|\widehat \lambda_k|}}\label{eq:gamma}.
\end{align}
Additionally, we consider $\widehat X^{(2)}$ such that
\begin{align}
    \|\widehat X^{(2)}\| \leq \epsilon.\label{eq:normX2}
\end{align}

We now discuss the upper bound of $\|\widehat X^{(2)}\|$.
From \eqref{eq:minamihata}, we have $|\underline{\widehat x}_{ij}^{(2)}| \leq u \cdot \mathrm{ufp}(\tau_j^{(1)})$.
Thus, we obtain
\begin{align}
    \|\widehat X^{(2)}\| \leq u\cdot\sqrt{\sum_{j=1}^n \left(\tau_j^{(1)}\right)^2}.\label{eq:xi}
\end{align}
If $\tau_{j}^{(1)} > 0$, from \eqref{eq:tau}, we have
\begin{align}
    \left(0.75 \cdot 2^{\beta}\right)^2 \sum_{j=1}^n 2^{2\lceil \log_2 w_{j}^{(1)} \rceil} &= \sum_{j=1}^n \left(\tau_j^{(1)}\right)^2.\label{eq:tau_ex}
\end{align}
Substituting \eqref{eq:tau_ex} into \eqref{eq:xi}, we obtain
\begin{align}
    \|\widehat X^{(2)}\|\leq  u \cdot \sqrt{\left(0.75 \cdot 2^{\beta}\right)^2 \sum_{j=1}^n 2^{2\lceil \log_2 w_{j}^{(1)} \rceil}}.
\end{align}
Thus, from \eqref{eq:normX2} and \eqref{eq:gamma}, we require
\begin{align}
    u \cdot \sqrt{(0.75 \cdot 2^{\beta})^2 \sum_{j=1}^n 2^{2\lceil \log_2 w_{j}^{(1)} \rceil}}\leq \sqrt{\delta\xi\cdot\frac{\min_{i \neq j}|\widehat \lambda_j - \widehat \lambda_i|}{\max_k|\widehat \lambda_k|}}.
\end{align}
Solving for $\beta$, we obtain
\begin{align}
    2^{\beta}\leq \frac{1}{0.75\cdot u}\cdot\sqrt{\frac{\delta\xi\cdot\min_{i\not=j}|\widehat\lambda_j-\widehat\lambda_i|}{\max_k|\widehat \lambda_k|\cdot\sum_{j=1}^n 2^{2\lceil \log_2 w_{j}^{(1)} \rceil}}}.\label{eq:beta2}
\end{align}
Then, we can set
\begin{align}
    {\beta}:= \left\lceil\log_2\left(\frac{1}{0.75\cdot u}\cdot\sqrt{\frac{\delta\xi\cdot\min_{i\not=j}|\widehat\lambda_j-\widehat\lambda_i|}{\max_k|\widehat \lambda_k|\cdot\sum_{j=1}^n 2^{2\lceil \log_2 w_{j}^{(1)} \rceil}}}\right)\right\rceil,\quad \alpha := \lceil-\log_2u+\log_2n-\beta\rceil\label{eq:proposed_beta}
\end{align}
for \eqref{eq:sigma} and \eqref{eq:tau}.

Using the computed values of $\alpha$ and $\beta$ together with \eqref{eq:a} and \eqref{eq:x}, we perform the slicing for $A$ and $\widehat{X}$ as follows:
\begin{align}
    A &= A^{(1)} + \dots + A^{(n_A-1)} + \underline{A}^{(n_A)}, 
    \quad \widehat{X} \leftarrow \widehat{X}^{(1)} .
\end{align}
After this procedure, when computing $A\widehat{X}$ with high accuracy, we have
\begin{align}
    A\widehat{X} 
    &= A^{(1)}\widehat{X} + \dots + A^{(n_A-1)}\widehat{X} + \underline{A}^{(n_A)}\widehat{X},\label{eq:mat_mul}
\end{align}
which can be implemented using $n_A$ matrix multiplications.  
In this case, it holds that
\[
    \mathrm{fl}(A^{(i)}\widehat{X}) = A^{(i)}\widehat{X}, \quad 1 \leq i \leq n_A - 1.
\]
The proposed method applies \eqref{eq:mat_mul} to the \textrm{accmul} in Algorithm 1.

\section{Numerical experiments}

All experiments were conducted on a computer with an Intel\textregistered{} Core\texttrademark{} i7-1270P CPU at 2.20 GHz and 32.0 GB of RAM, running Windows 11 Pro. MATLAB R2024a was used.
The methods were implemented in MATLAB without the MEX (MATLAB Executable, using C or FORTRAN) function.

In this section, we verify the efficiency of the proposed method through numerical experiments.
For comparison, we consider the following cases:
\begin{itemize}
    \item Previous study: $n_A = n_X = 2,3,4$,
    \item Proposed method: $n_X=1, n_A = 2, \dots, 6$.
\end{itemize}
In the previous study, since $n_A = n_X = 3$, the accurate matrix multiplication requires six matrix multiplications.  
Therefore, the number of matrix multiplications required for iterative refinement is $8 \times$ (the number of iterations).  
For the proposed method, the number of matrix multiplications required for iterative refinement is $(n_A + 2) \times$ (the number of iterations).  
Thus, when $n_A = 6$, the computational cost of the proposed method is equivalent to that of the previous study.

The test matrix was generated using the following MATLAB code.
\begin{verbatim}
    A = gallery('randsvd',n,-cnd,3,n-1,n-1,1);
\end{verbatim}
The approximate eigenvalue decomposition was computed using the \texttt{eig} function.

\subsection{Implementation}

In this section, we describe the implementation method of the proposed approach.
Mathematically, the determination of $\alpha$ and $\beta$ follows a defined procedure.
However, in many cases, the error estimation tends to be overly conservative, meaning that the values may not be optimal.

Here, we set
\begin{align}
     \alpha = -(\log_2 u + \beta) + \lceil \log_2 \sqrt{n} \rceil,\quad 
     \beta =\left\lfloor\log_2\left(\frac{1}{0.75\cdot u}\cdot\sqrt{\frac{n\delta\cdot\min_{i\not=j}|\widehat\lambda_j-\widehat\lambda_i|}{\max_k|\widehat \lambda_k|\cdot\sum_{j=1}^n 2^{2\lceil \log_2 w_{j}^{(1)} \rceil}}}\right)\right\rfloor
\end{align}
from~\eqref{eq:proposed_beta}.
These values do not strictly satisfy the condition \eqref{eq:alp_beta}.
As a result, for $1 \leq i < n_A, 1 \leq j < n_X$, there is a possibility that $A^{(i)}\widehat X^{(j)} \neq \mathrm{fl}(A^{(i)}\widehat X^{(j)})$.
On the other hand, when $A^{(i)}\widehat X^{(j)} = \mathrm{fl}(A^{(i)}\widehat X^{(j)})$, this approach improves computational accuracy.
In particular, as the matrix size $n$ increases, this method helps to suppress excessive error overestimation.

Furthermore, the value of $\beta$ is set relatively large to account for the overestimation of the error evaluations.
These parameters can be adjusted depending on the specific problem being handled.

\subsection{Test for random matrices}

In this section, we present the results of numerical experiments on the forward error $\|X - \widehat X\|$ of the eigenvector matrix $\widehat X$ obtained using the proposed method.
\begin{figure}[htbp]
    \centering
    \begin{subfigure}{0.32\textwidth}
        \centering
        \includegraphics[width=\textwidth]{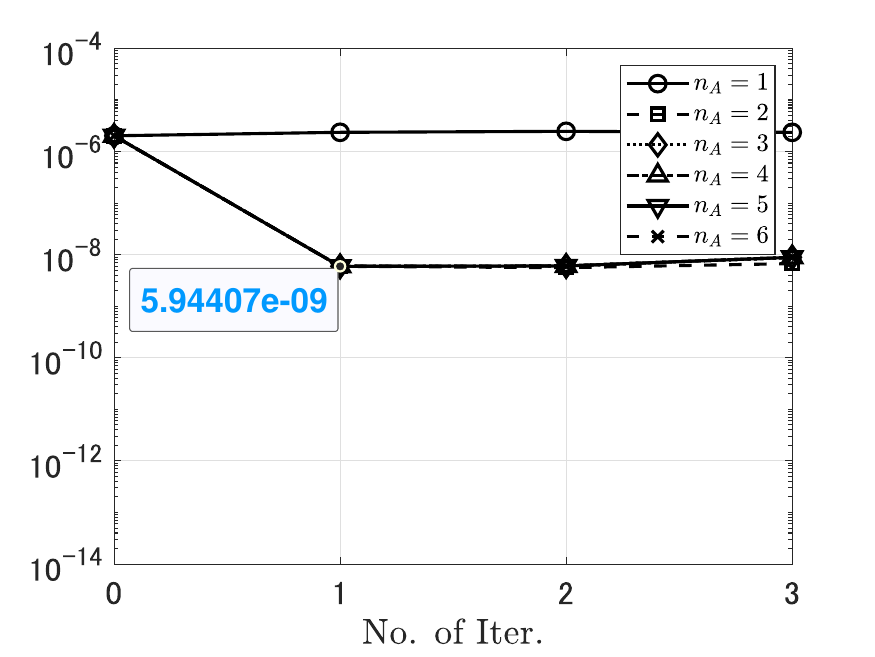} 
        \subcaption{$n=1024, \delta=10^{-8}$}\label{his:a}
    \end{subfigure}
    \hfill
    \begin{subfigure}{0.32\textwidth}
        \centering
        \includegraphics[width=\textwidth]{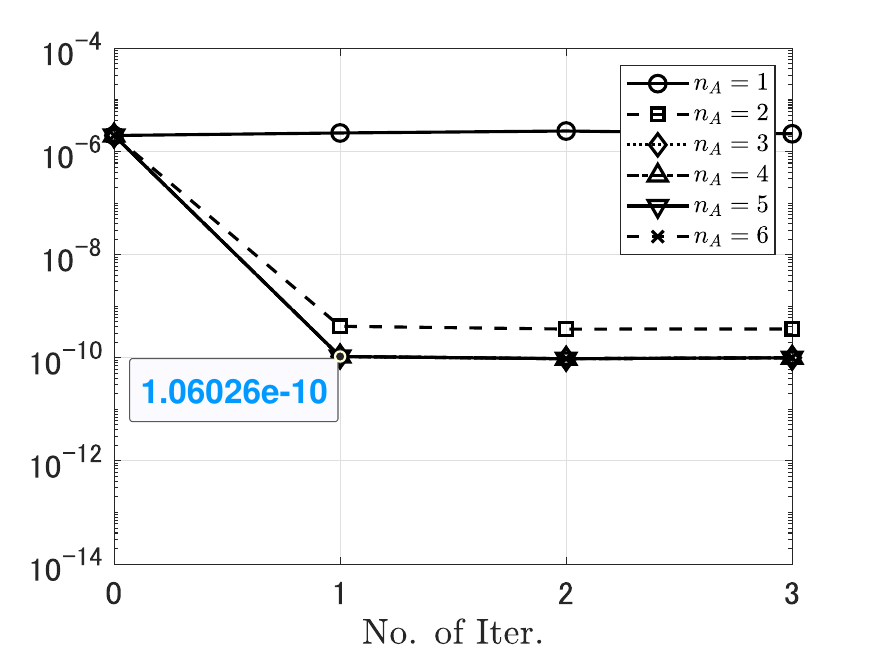} 
        \subcaption{$n=1024, \delta=10^{-10}$}\label{his:b}
    \end{subfigure}
    \hfill
    \begin{subfigure}{0.32\textwidth}
        \centering
        \includegraphics[width=\textwidth]{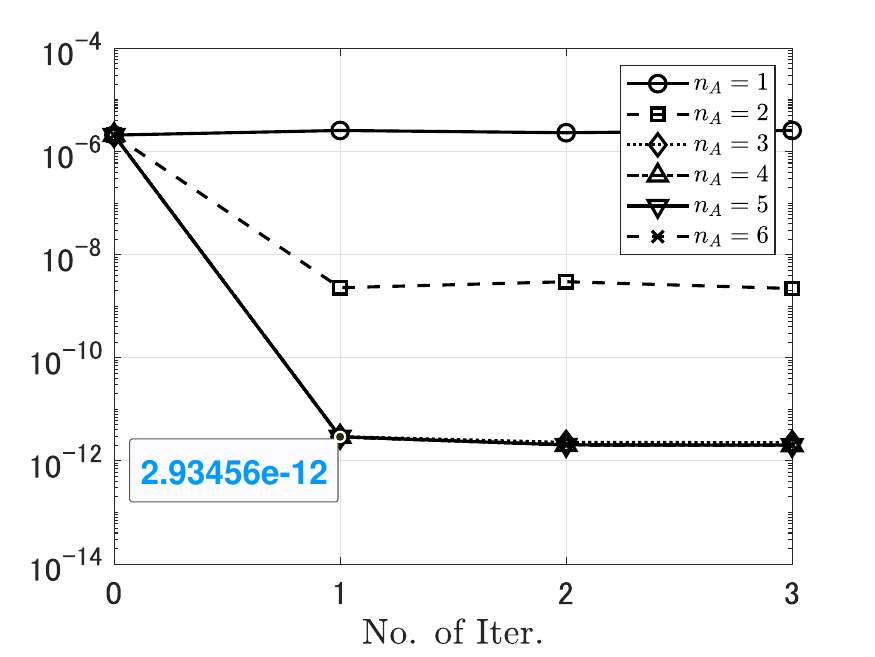} 
        \subcaption{$n=1024, \delta=10^{-12}$}\label{his:c}
    \end{subfigure}\\

    \begin{subfigure}{0.32\textwidth}
        \centering
        \includegraphics[width=\textwidth]{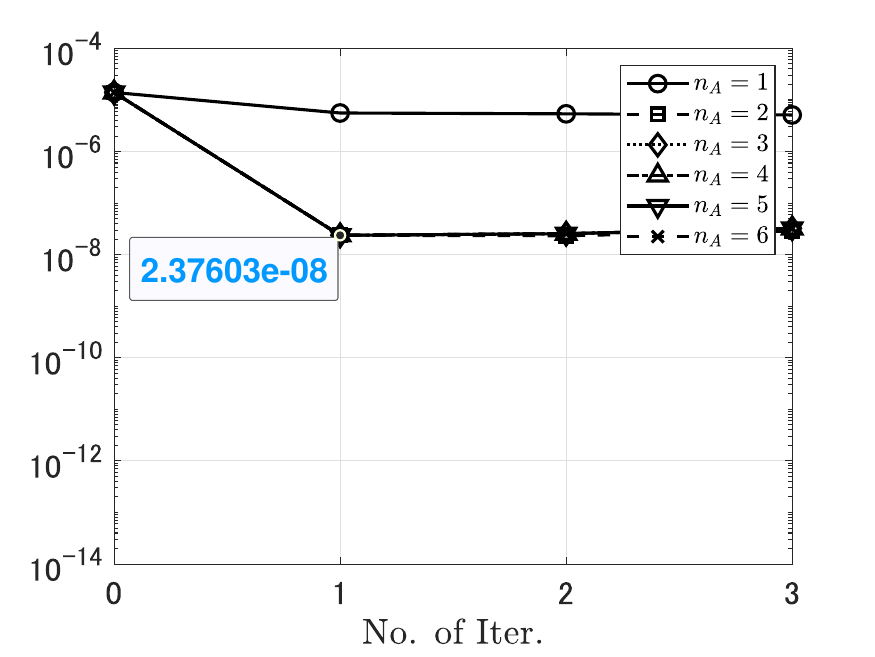} 
        \subcaption{$n=4096, \delta=10^{-8}$}\label{his:d}
    \end{subfigure}
    \hfill
    \begin{subfigure}{0.32\textwidth}
        \centering
        \includegraphics[width=\textwidth]{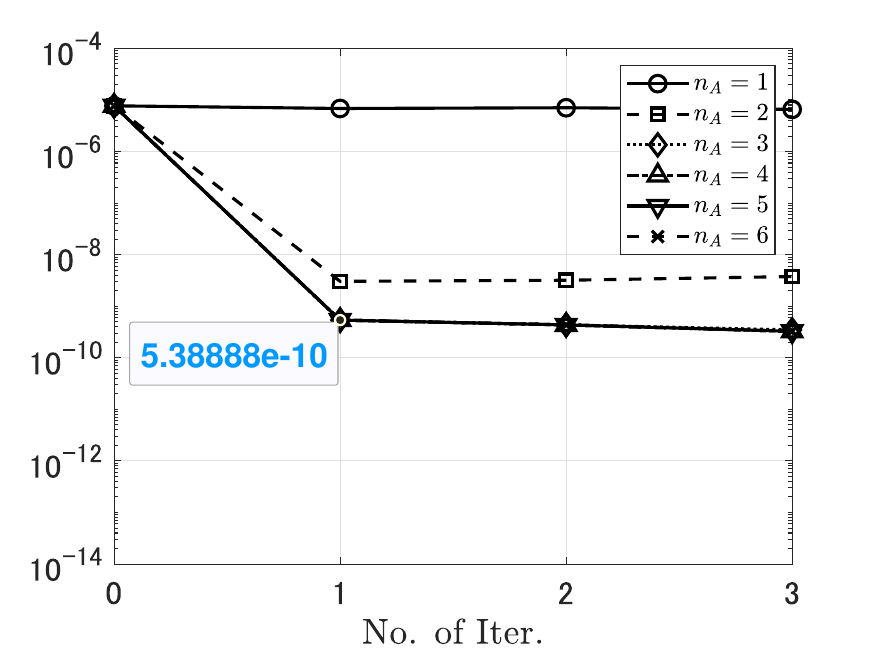} 
        \subcaption{$n=4096, \delta=10^{-10}$}\label{his:e}
    \end{subfigure}
    \hfill
    \begin{subfigure}{0.32\textwidth}
        \centering
        \includegraphics[width=\textwidth]{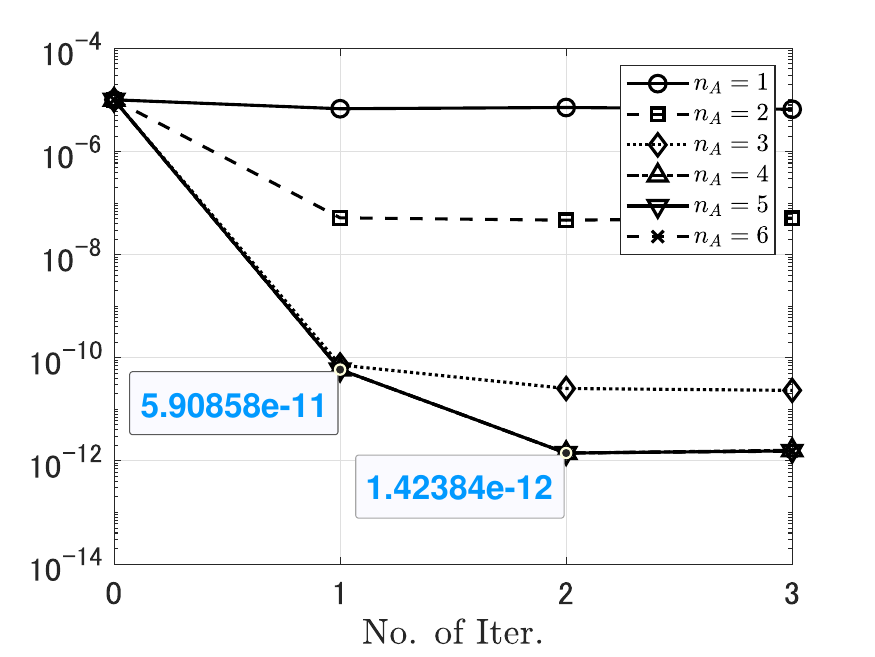} 
        \subcaption{$n=4096, \delta=10^{-12}$}\label{his:f}
    \end{subfigure}\\

    \begin{subfigure}{0.32\textwidth}
        \centering
        \includegraphics[width=\textwidth]{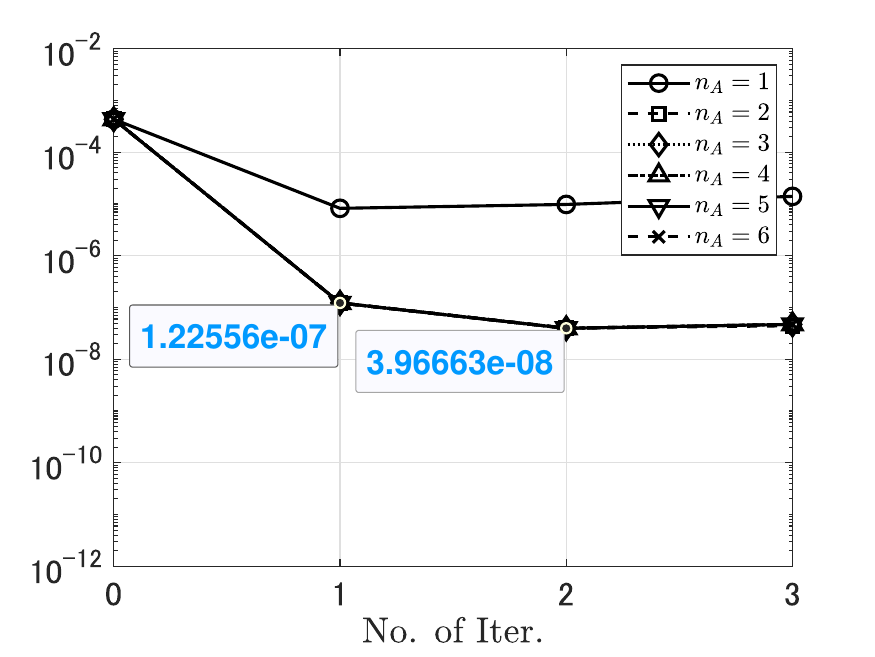} 
        \subcaption{$n=8192, \delta=10^{-8}$}\label{his:h}
    \end{subfigure}
    \hfill
    \begin{subfigure}{0.32\textwidth}
        \centering
        \includegraphics[width=\textwidth]{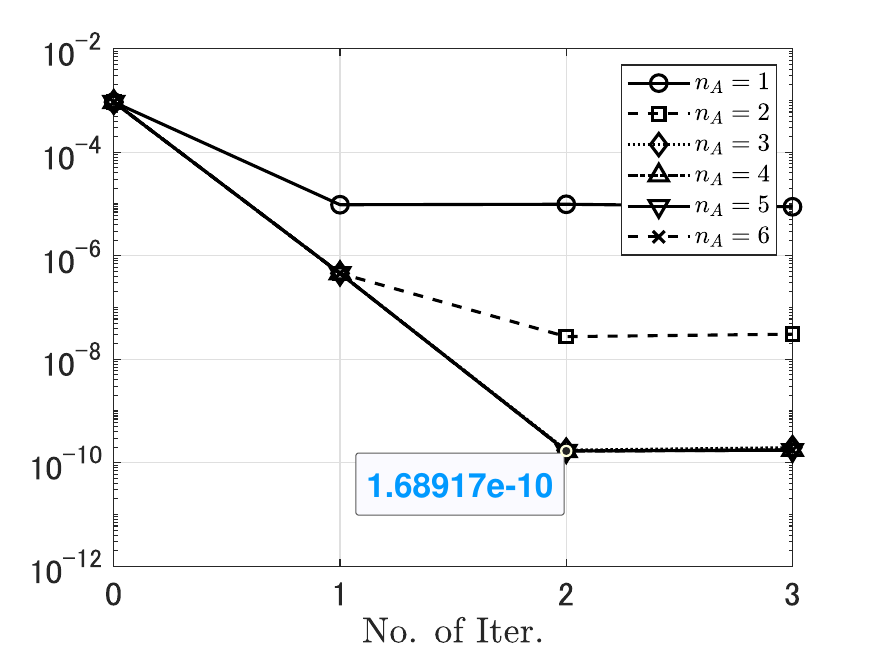} 
        \subcaption{$n=8198, \delta=10^{-10}$}\label{his:i}
    \end{subfigure}
    \hfill
    \begin{subfigure}{0.32\textwidth}
        \centering
        \includegraphics[width=\textwidth]{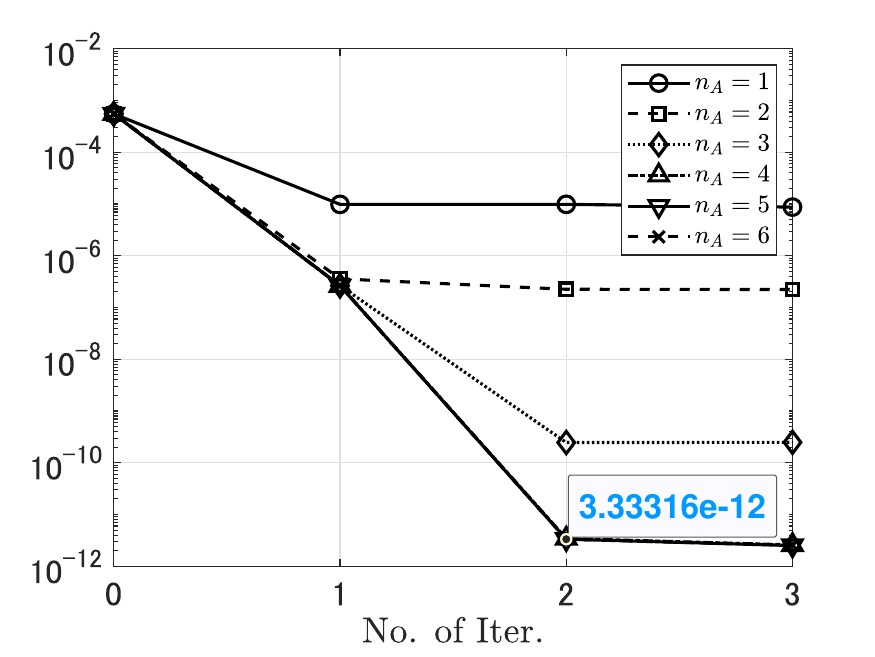} 
        \subcaption{$n=8198, \delta=10^{-12}$}\label{his:j}
    \end{subfigure}
    \caption{Convergence history of forward error $\|X - \widehat X\|$ for proposed method ($\mathtt{cnd}=10^{10}$).}\label{fig:his}
\end{figure}

Figure~\ref{fig:his} shows the convergence history of the proposed method.
We consider nine cases by combining $n=1024,4096,8192$ with $\delta=10^{-8},10^{-10},10^{-12}$.
For the iteration count $0$, the initial accuracy (computed by MATLAB's \texttt{eig} function) is shown.

From Figure~\ref{fig:his}, we observe that as the matrix size increases or the problem becomes ill-conditioned, both $n_A$ and the number of iterations increase.
However, a value of $n_A$ between 2 and 4 is sufficient to improve the accuracy of the eigenvector matrix to a satisfactory level.
Additionally, in most cases, one or two iterations were sufficient to achieve the expected accuracy.
These results demonstrate that the proposed method reduces the computational cost by a factor of 4/3 to 2 compared to the method in~\cite{uchino2024high}.
Furthermore, from Figure~\ref{fig:his}, we observe that
\begin{align}
    \|X - \widehat X\| = \mathcal{O}(\delta),
\end{align}
indicating that the proposed method successfully computes approximate eigenvectors with the specified forward error.

\begin{figure}[htbp]
    \centering
    \begin{subfigure}{0.32\textwidth}
        \centering
        \includegraphics[width=\textwidth]{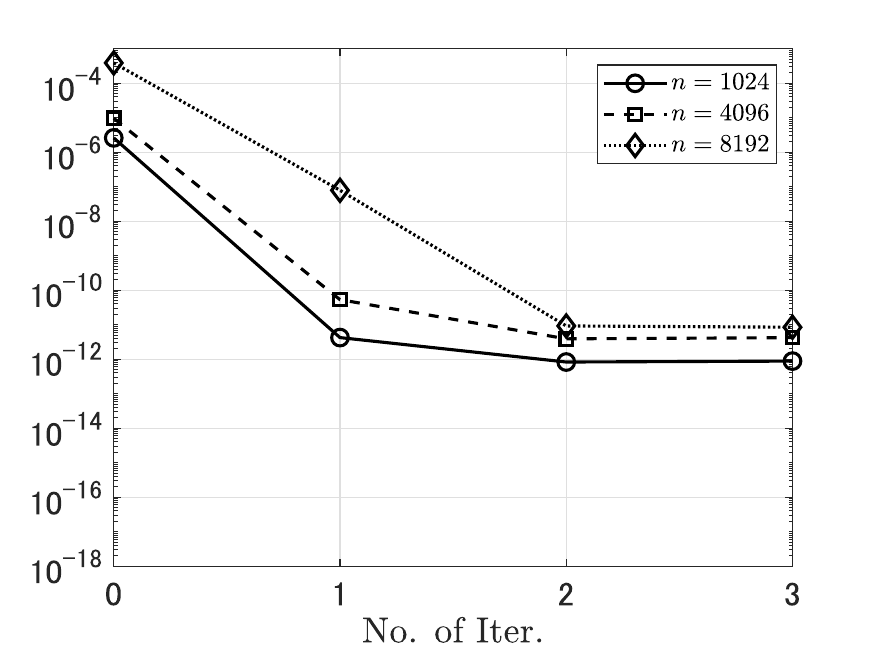} 
        \subcaption{$n_A=n_X=2$}\label{his-u:a}
    \end{subfigure}
    \hfill
    \begin{subfigure}{0.32\textwidth}
        \centering
        \includegraphics[width=\textwidth]{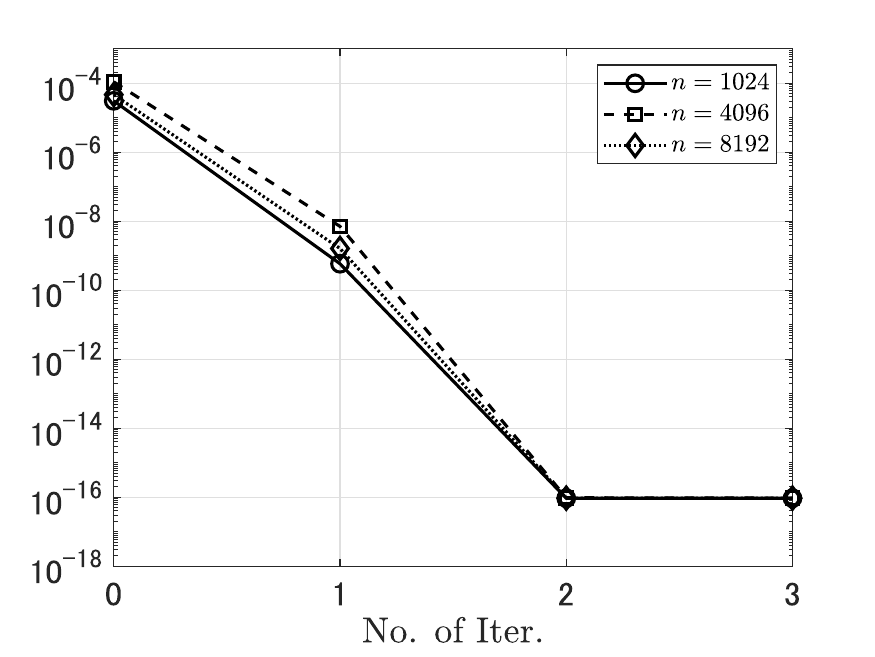} 
        \subcaption{$n_A=n_X=3$}\label{his-u:b}
    \end{subfigure}
    \hfill
    \begin{subfigure}{0.32\textwidth}
        \centering
        \includegraphics[width=\textwidth]{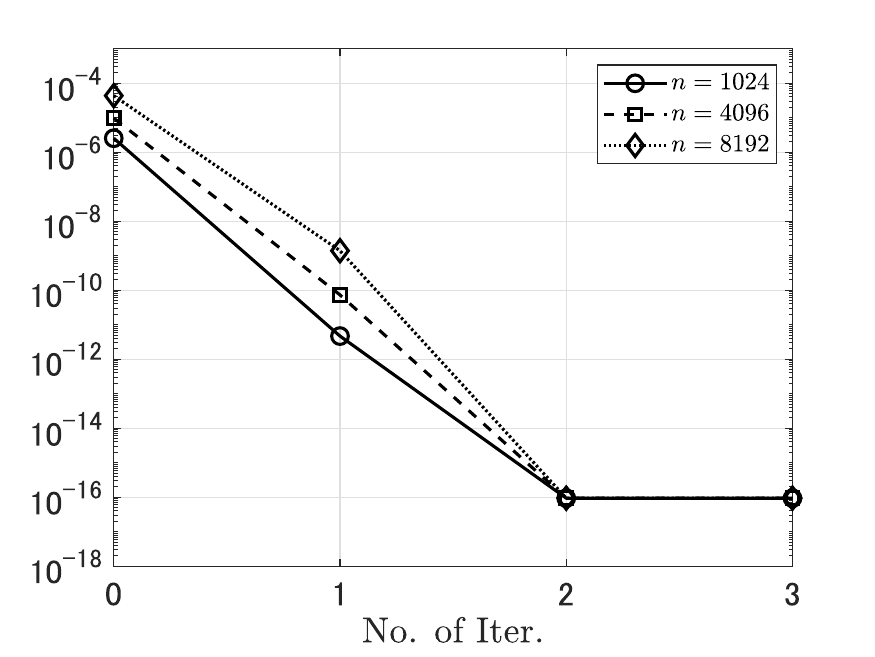} 
        \subcaption{$n_A=n_X=4$}\label{his-u:c}
    \end{subfigure}
    \caption{Convergence history of forward error $\|X - \widehat X\|$ for method in~\cite{uchino2024high} ($\mathtt{cnd}=10^{10}$).}
    \label{fig_uchino}
\end{figure}

Next, we evaluate the computational accuracy of the previous study~\cite{uchino2024high}.
Figure~\ref{fig_uchino} shows the forward error $\|X - \widehat X\|$ for $n_A = n_X = 3$.
Comparing Figures~\ref{fig:his} and \ref{fig_uchino}, we find that the convergence slowdown due to the splitting of $\widehat X$ in the proposed method is minimal.
This confirms that the proposed method efficiently computes approximate solutions with a specified forward error.

\subsection{Sparse matrices}

Here, we conduct experiments using real-world matrices from the ELSES matrix library. These matrices are generated from electronic structure calculations. For details, please refer to \cite{hoshi2006large,hoshi2019numerical}.
The name of the matrix contains the size. For example, the name ``VCNT800std'' contains matrices of size $n=800$.

We set
\begin{align}
     \alpha &= -(\log_2 u + \beta) + \lceil \log_2 k \rceil,\\ 
     \beta  &=\left\lfloor\log_2\left(\frac{\gamma \cdot \min(k,\sqrt{n})}{0.75\cdot u}\cdot\sqrt{\frac{1}{\sum_{j=1}^n 2^{2\lceil \log_2 w_{j}^{(1)} \rceil}}}\right)\right\rfloor,\\ 
     \gamma &=\sqrt{\delta\cdot\frac{\min_{i \neq j}|\widehat \lambda_j - \widehat \lambda_i|}{\max_k|\widehat \lambda_k|}},
\end{align}
where $k$ is the maximum number of nonzero elements in any row of matrix $A$.

\begin{table}[htbp]
 \centering
   \caption{Test for sparse matrices of the ELSES matrix library}\label{tab:time2}
  \begin{tabular}{rcrrr|rr|rrr}
     & & & \multicolumn{7}{c}{$\|\widehat X-X\|$} \\\cmidrule(lr){4-10}
     & & & \multicolumn{2}{c}{$\delta=10^{-6}$} & \multicolumn{2}{c}{$\delta=10^{-8}$} & \multicolumn{3}{c}{$\delta=10^{-10}$} \\\cmidrule(lr){4-5}\cmidrule(lr){6-7}\cmidrule(lr){8-10}
    \multicolumn{1}{c}{name} & $n_A$ & $\|X_0-X\|$ & Iter. = 1 & Iter. = 2 & Iter. = 1 & Iter. = 2 & Iter. = 1 & Iter. = 2 & Iter. = 3 \\\midrule
    \multirow{1}{*}{VCNT800std} & $2$ & \multirow{1}{*}{1.71e-03} &  6.05e-06 & -- & 1.25e-06 & 5.12e-08 & 1.25e-06 & 3.11e-10 & --\\
    \multirow{1}{*}{VCNT2000std} & $2$ & \multirow{1}{*}{2.86e-03} & 5.12e-06 & -- & 5.04e-06 & 5.01e-08 & 5.04e-06 & 6.06e-10 & --\\
    \multirow{1}{*}{VCNT4000std} & $2$ & \multirow{1}{*}{2.68e-02} & 2.65e-04 & 4.70e-06 & 2.65e-04 & 6.41e-08 & 2.65e-04 & 6.49e-08 & 4.14e-10\\
  \end{tabular}
\end{table}

\begin{table}[htbp]
 \centering
   \caption{Test for SuiteSparse matrix collection with $1000\leq n\leq 5000, \delta =10^{-6}$, and $n_A=2$}\label{tab:sparse}
  \begin{tabular}{rcc|cccccc}
      &  & & \multicolumn{6}{c}{$\|\widehat X-X\|$} \\ \cmidrule(lr){4-9}
     id & $n$ & $\|X_0-X\|$ & Iter. = 1 & Iter. = 2 & Iter. = 3 & Iter. = 4 & Iter. = 5 & Iter. = 6\\ \midrule
182 & 1138 & 1.67e-03 & 1.12e-06 & 7.82e-08 & 9.40e-08 \\ 
183 & 1141 & 1.23e-03 & 8.43e-07 & 1.83e-07 \\ 
483 & 2000 & 2.12e-02 & 1.80e-04 & 6.53e-07 & 7.33e-07 \\ 
484 & 2000 & 3.64e-03 & 6.09e-06 & 3.79e-07 & 8.22e-07 \\ 
485 & 2000 & 3.26e-03 & 4.34e-06 & 8.75e-07 \\ 
486 & 2000 & 5.58e-03 & 1.24e-05 & 4.99e-07 & 5.30e-07 \\ 
487 & 2000 & 3.67e-03 & 6.87e-06 & 1.47e-06 \\ 
488 & 2000 & 8.83e-04 & 1.61e-06 & 4.84e-06 \\ 
489 & 2000 & 1.32e-03 & 1.81e-06 & 2.24e-06 \\ 
490 & 2000 & 1.62e-03 & 2.43e-06 & 2.46e-06 \\ 
492 & 2000 & 4.83e-03 & 9.14e-06 & 4.41e-07 & 1.15e-06 \\ 
493 & 2000 & 1.55e-03 & 2.80e-06 & 3.21e-06 \\ 
494 & 2000 & 8.17e-03 & 1.46e-05 & 1.04e-07 & 1.37e-07 \\ 
495 & 2000 & 2.81e-03 & 5.50e-06 & 6.20e-08 & 1.85e-08 \\ 
496 & 2000 & 5.09e-03 & 1.16e-05 & 9.69e-08 & 8.53e-08 \\ 
497 & 2000 & 3.19e-01 & 5.28e-02 & 3.46e-03 & 2.22e-05 & 3.93e-07 & 8.06e-07 \\ 
498 & 2000 & 2.66e-02 & 4.11e-04 & 2.48e-06 & 9.19e-07 \\ 
499 & 2000 & 1.68e-02 & 1.13e-04 & 1.95e-06 & 2.40e-06 \\ 
500 & 2000 & 7.55e-03 & 2.37e-05 & 1.51e-06 & 1.40e-06 \\ 
501 & 2000 & 3.65e-03 & 6.19e-06 & 1.45e-06 \\ 
503 & 2000 & 4.86e-03 & 9.58e-06 & 1.56e-06 \\ 
504 & 2000 & 8.13e-03 & 2.64e-05 & 1.72e-06 & 1.30e-06 \\ 
505 & 2000 & 8.21e-03 & 2.68e-05 & 2.02e-06 & 2.19e-06 \\ 
521 & 5000 & 1.03e-01 & 4.28e-03 & 1.80e-05 & 7.79e-09 & 3.48e-08 \\ 
522 & 5000 & 4.17e-01 & 1.09e-01 & 1.74e-02 & 4.22e-04 & 2.14e-06 & 1.66e-06 \\ 
523 & 5000 & 2.83e-01 & 3.84e-02 & 1.70e-03 & 3.90e-06 & 4.11e-06 \\ 
758 & 2146 & 6.13e-01 & 2.50e-01 & 1.08e-01 & 7.74e-03 & 2.91e-05 & 2.05e-06 & 2.08e-06 \\ 
808 & 4350 & 2.10e-01 & 2.41e-02 & 4.68e-04 & 6.08e-07 & 6.22e-07 \\ 
1442 & 1357 & 8.85e-03 & 1.94e-05 & 2.19e-07 & 1.60e-07 \\ 
2211 & 2000 & 6.62e-03 & 1.66e-05 & 7.67e-08 & 7.52e-08 \\ 
2409 & 4720 & 1.91e-02 & 1.67e-04 & 6.87e-08 & 5.15e-08 \\ 
2411 & 4253 & 5.62e-02 & 1.42e-03 & 1.98e-06 & 6.71e-08 & 6.28e-08 \\ 
2423 & 1866 & 4.25e-01 & 1.18e-01 & 2.13e-02 & 6.11e-04 & 1.46e-07 & 2.56e-08 \\ 
2465 & 1024 & 1.80e-03 & 1.44e-06 & 1.78e-07 \\ 
2466 & 2048 & 3.21e-03 & 3.80e-06 & 2.01e-07 & 1.54e-07 \\ 
2467 & 4096 & 1.89e-02 & 1.51e-04 & 2.21e-07 & 1.76e-07 \\ 
2663 & 2427 & 1.45e-01 & 8.76e-03 & 8.62e-05 & 3.03e-07 & 2.45e-07 \\ 
2860 & 1000 & 1.96e-02 & 1.90e-04 & 7.23e-08 & 9.54e-08 \\ 
2862 & 2000 & 5.60e-03 & 3.52e-05 & 1.80e-06 & 1.39e-06 \\ 
2876 & 1077 & 5.48e-01 & 3.43e-01 & 1.39e-01 & 1.17e-02 & 1.32e-04 & 2.40e-07 & 1.82e-07 \\ 
2881 & 1600 & 7.12e-02 & 2.57e-03 & 9.38e-06 & 4.37e-07 & 5.90e-07 \\ 
2882 & 1600 & 5.82e-01 & 2.77e-01 & 1.54e-01 & 1.66e-02 & 1.61e-04 & 5.94e-07 & 3.68e-07 \\ 
2883 & 1599 & 2.29e-02 & 1.98e-04 & 1.31e-06 & 2.20e-06 \\ 
2884 & 1593 & 2.96e-03 & 7.94e-06 & 2.73e-06 \\ 
2891 & 2414 & 7.59e-02 & 3.02e-03 & 7.18e-06 & 2.33e-06 \\ 
\end{tabular}
\end{table}

\begin{table}[htbp]
 \centering
   \caption{Test for SuiteSparse matrix collection with $1000\leq n\leq 5000, \delta =10^{-10}$, and $n_A=2$}\label{tab:sparse3}
  \begin{tabular}{rcc|ccccccc}
      &  & & \multicolumn{7}{c}{$\|\widehat X-X\|$} \\ \cmidrule(lr){4-10}
     id & $n$ & $\|X_0-X\|$ & Iter. = 1 & Iter. = 2 & Iter. = 3 & Iter. = 4 & Iter. = 5 & Iter. = 6 & Iter. = 7\\ \midrule
182 & 1138 & 1.67e-03 & 1.12e-06 & 1.78e-11 & 2.96e-11 \\ 
183 & 1141 & 1.23e-03 & 8.48e-07 & 1.44e-11 & 1.78e-11 \\ 
483 & 2000 & 2.12e-02 & 1.80e-04 & 3.77e-08 & 1.91e-10 & 1.91e-10 \\ 
484 & 2000 & 3.64e-03 & 6.09e-06 & 1.12e-10 & 9.00e-11 \\ 
485 & 2000 & 3.26e-03 & 4.33e-06 & 1.21e-10 & 1.03e-10 \\ 
486 & 2000 & 5.58e-03 & 1.24e-05 & 1.79e-10 & 1.03e-10 \\ 
487 & 2000 & 3.67e-03 & 6.87e-06 & 1.13e-10 & 1.26e-10 \\ 
488 & 2000 & 8.83e-04 & 5.32e-07 & 2.52e-10 & 3.73e-10 \\ 
489 & 2000 & 1.32e-03 & 9.54e-07 & 1.54e-10 & 1.73e-10 \\ 
490 & 2000 & 1.62e-03 & 1.32e-06 & 1.48e-10 & 9.88e-11 \\ 
492 & 2000 & 4.83e-03 & 9.13e-06 & 4.69e-10 & 3.12e-10 \\ 
493 & 2000 & 1.55e-03 & 1.02e-06 & 4.88e-10 & 2.37e-10 \\ 
494 & 2000 & 8.17e-03 & 1.46e-05 & 1.59e-10 & 1.19e-11 & 1.19e-11 \\ 
495 & 2000 & 2.81e-03 & 5.53e-06 & 2.75e-11 & 1.30e-11 \\ 
496 & 2000 & 5.09e-03 & 1.16e-05 & 1.77e-10 & 3.83e-12 & 3.83e-12 \\ 
497 & 2000 & 3.19e-01 & 5.28e-02 & 3.45e-03 & 2.15e-05 & 1.92e-10 & 1.00e-10 \\ 
498 & 2000 & 2.66e-02 & 4.11e-04 & 1.26e-07 & 1.01e-10 & 9.35e-11 \\ 
499 & 2000 & 1.68e-02 & 1.13e-04 & 1.46e-08 & 6.58e-10 & 6.58e-10 \\ 
500 & 2000 & 7.55e-03 & 2.37e-05 & 6.30e-10 & 3.11e-10 \\ 
501 & 2000 & 3.65e-03 & 6.20e-06 & 3.39e-10 & 3.80e-10 \\ 
503 & 2000 & 4.86e-03 & 9.58e-06 & 1.18e-10 & 1.17e-10 \\ 
504 & 2000 & 8.13e-03 & 2.64e-05 & 7.76e-10 & 2.41e-10 \\ 
505 & 2000 & 8.21e-03 & 2.68e-05 & 8.70e-10 & 8.61e-10 \\ 
521 & 5000 & 1.03e-01 & 4.28e-03 & 1.80e-05 & 2.90e-10 & 5.22e-11 \\ 
522 & 5000 & 4.17e-01 & 1.09e-01 & 1.74e-02 & 4.21e-04 & 7.00e-08 & 2.82e-10 & 2.76e-10 \\ 
523 & 5000 & 2.83e-01 & 3.84e-02 & 1.70e-03 & 5.38e-06 & 1.77e-10 & 1.77e-10 \\ 
758 & 2146 & 6.13e-01 & 2.50e-01 & 1.08e-01 & 7.74e-03 & 2.93e-05 & 1.33e-09 & 7.15e-11 & 7.39e-11 \\ 
808 & 4350 & 2.10e-01 & 2.41e-02 & 4.68e-04 & 2.96e-07 & 3.53e-11 & 3.53e-11 \\ 
1442 & 1357 & 8.85e-03 & 1.93e-05 & 3.74e-10 & 2.06e-11 & 2.06e-11 \\ 
2211 & 2000 & 6.62e-03 & 1.66e-05 & 3.24e-10 & 4.00e-11 \\ 
2411 & 4253 & 5.62e-02 & 1.42e-03 & 1.97e-06 & 1.10e-11 & 1.10e-11 \\ 
2409 & 4720 & 1.91e-02 & 1.67e-04 & 2.39e-08 & 1.87e-11 & 1.61e-11 \\ 
2423 & 1866 & 4.25e-01 & 1.18e-01 & 2.13e-02 & 6.11e-04 & 1.48e-07 & 3.37e-10 & 3.37e-10 \\ 
2465 & 1024 & 1.80e-03 & 1.44e-06 & 8.32e-11 & 4.72e-11 \\ 
2466 & 2048 & 3.21e-03 & 3.80e-06 & 1.12e-10 & 4.76e-11 \\ 
2467 & 4096 & 1.89e-02 & 1.51e-04 & 2.29e-08 & 4.79e-11 & 4.79e-11 \\ 
2663 & 2427 & 1.45e-01 & 8.76e-03 & 8.62e-05 & 8.41e-09 & 1.22e-10 & 1.22e-10 \\ 
2860 & 1000 & 1.96e-02 & 1.90e-04 & 5.31e-08 & 1.59e-11 & 1.59e-11 \\ 
2862 & 2000 & 5.60e-03 & 3.57e-05 & 9.36e-10 & 1.07e-10 \\ 
2876 & 1077 & 5.48e-01 & 3.43e-01 & 1.39e-01 & 1.17e-02 & 1.33e-04 & 1.34e-08 & 1.44e-11 & 1.56e-11 \\ 
2881 & 1600 & 7.12e-02 & 2.57e-03 & 1.04e-05 & 3.19e-10 & 1.08e-11 & 1.08e-11 \\ 
2882 & 1600 & 5.82e-01 & 2.77e-01 & 1.54e-01 & 1.66e-02 & 1.61e-04 & 3.03e-08 & 1.32e-10 & 1.41e-10 \\ 
2883 & 1599 & 2.29e-02 & 1.98e-04 & 4.43e-08 & 1.15e-10 & 1.06e-10 \\ 
2884 & 1593 & 2.96e-03 & 7.73e-06 & 4.51e-11 & 6.70e-11 \\ 
2891 & 2414 & 7.59e-02 & 3.02e-03 & 7.18e-06 & 1.22e-10 & 1.22e-10 \\ 
  \end{tabular}
\end{table}

From these results, if the given matrix is a sparse matrix arising from practical applications, a two-way partition of $A$ is sufficient.
Thus, the computation required for high-precision matrix multiplication reduces to two standard matrix multiplications, enabling acceleration.
The results also indicate that, for many problems, convergence to a prescribed forward error can be achieved.
Consequently, for practical problems, it becomes possible to compute eigenvectors that attain a specified forward error while keeping the cost of high-precision arithmetic low.
Therefore, this approach is expected to be applicable to future high-precision and mixed-precision computations.

\section{Summary}

In this paper, we have proposed an iterative refinement method for real symmetric eigenvalue problems, focusing on controlling forward errors of eigenvectors.  
Unlike conventional approaches that primarily minimize backward errors, our method is targeted at a user-specified forward error tolerance.  
By adapting the Ogita--Aishima iterative refinement framework and incorporating the Ozaki scheme for cost-effective high-precision matrix multiplications, we reduce computational cost while preserving accuracy.  
Numerical experiments on both random and real-world sparse matrices demonstrate that the proposed method attains the target forward error with fewer matrix multiplications than previous approaches.  
These results confirm that the method realized an efficient and reliable eigensolver for applications requiring forward error guarantees.  
Although clustered eigenvalues are not addressed in this paper, a natural extension of iterative refinement to handle such cases would allow our framework to cover them.  
Furthermore, the proposed method can be extended in a straightforward manner to Hermitian matrices.






\backmatter

\bmhead{Acknowledgements}

This work was partially supported by JSPS KAKENHI Grant Numbers 23H03410 and 25K03126. 

\bibliography{mybib}

\end{document}